\documentclass[11pt, reqno]{amsart}

\usepackage[T1]{fontenc}
\usepackage{amsmath}
\usepackage{amssymb}
\usepackage{amsthm}
\usepackage[left=3.5cm, right=3.5cm, paperheight=11.9in]{geometry}
\usepackage{hyperref}
\usepackage{fancyhdr}
\usepackage{enumitem}
\usepackage{comment}
\usepackage{nicefrac}
\usepackage{bm}
\usepackage{mathrsfs}
\usepackage{graphicx}
\usepackage[utf8]{inputenc}
\usepackage{cancel}
\usepackage{mathtools}

\author[P.Y. Bienvenu]{Pierre-Yves Bienvenu}
\address{Johann Radon Institute for Computational and Applied Mathematics, Austrian Academy of Sciences, Linz, Austria}
\email{pierre.bienvenu@oeaw.ac.at}

\author[P.~Leonetti]{Paolo Leonetti}
\address{
Department of Economics, 
Universit\`a degli Studi dell'Insubria, via Monte Generoso 71, 21100 Varese, Italy}
\email{leonetti.paolo@gmail.com}

\newtheorem{thm}{Theorem}[section]
\newtheorem{cor}[thm]{Corollary}
\newtheorem{lem}[thm]{Lemma}
\newtheorem{lemma}[thm]{Lemma}

\newtheorem{prop}[thm]{Proposition}

\theoremstyle{definition} 

\newtheorem{rmk}[thm]{Remark}
\newtheorem{question}[thm]{Question}

\numberwithin{equation}{section}

\newcommand{\Z}{\mathbb{Z}}

\newcommand{\cR}{\mathcal{R}}
\newcommand{\N}{\mathbb{N}}
\newcommand{\I}{\mathcal{I}}
\newcommand{\Pn}{\mathcal{P}(\N)}
\newcommand{\Fin}{\mathrm{Fin}}

\newcommand{\sd}{\mathbin{\triangle}}

\newcounter{smallromans}

{\end{list}}

\AtBeginDocument{%
   \def\MR#1{}
}

\hypersetup{
    pdftitle={Sumsets modulo ideals: Baire category and descriptive complexity},
    pdfauthor={Paolo Leonetti and Pierre Bienvenu},
    pdfmenubar=false,
    pdffitwindow=true,
    pdfstartview=FitH,
    colorlinks=true,
    linkcolor=blue,
    citecolor=green,
    urlcolor=cyan
}

\providecommand{\MR}[1]{}

\providecommand{\MR}{\relax\ifhmode\unskip\space\fi MR }

\begin{document}
\title{Sumsets modulo ideals: Baire category and descriptive complexity}



\begin{abstract} 
Given an ideal $\mathcal I$ on the nonnegative integers, we study the
family $\Sigma(\mathcal I)$ of all sets for which every $\mathcal{I}$-small modification is not a nontrivial sumset $A+B$. 
We prove that $\Sigma(\mathcal I)$ is comeager whenever $\mathcal{I}$ is an Erd\H{o}s--Ulam ideal generated by a nonincreasing sequence of positive weights. This answers an open question for the ideal $\mathcal{Z}$ of asymptotic density zero sets in [Proc. Amer. Math. Soc. \textbf{151} (2023), 3681--3689]. 

We also show that $\Sigma(\mathcal I)$ is either empty or a
dense $G_\delta$-complete set whenever $\mathcal{I}$ is an $F_\sigma$ ideal which is invariant under finite modifications. This implies that the family of nontrivial sumsets $A+B$ is $F_\sigma$-complete (hence, Borel). On the other hand, we construct a meager ideal $\mathcal{I}$ for which $\Sigma(\mathcal I)$ is neither meager nor comeager. Finally, we
prove that $\Sigma(\mathcal Z)$ is coanalytic and not
$F_{\sigma\delta}$.
\end{abstract}
\maketitle
\thispagestyle{empty}

\section{Introduction}
\subsection{Sumsets and irreducibility}
A subset $A$ of the nonnegative integers $\mathbb{N}$ is said to be \emph{irreducible}
if it is not a nontrivial sumset, i.e., if there are no $B,C\subseteq \mathbb{N}$ such that $|B|,|C|\ge 2$ and $A=B+C:=\{b+c: b \in B, c \in C\}$. 
Moreover, $A\subseteq \mathbb{N}$ is called \emph{totally irreducible} if there are no sets $B,C\subseteq \mathbb{N}$ such that $|B|,|C|\ge 2$ and $A\bigtriangleup (B+C)$ is finite, cf. \cite{MR4689545, MR58655}. Let $\Sigma$ and $\Sigma_\star$ be the families of irreducible and totally irreducible sets, respectively, so that $\Sigma_\star\subseteq \Sigma$.
Let $\cR$ and $\cR_*$ respectively be their complements to $\Pn$.
Wirsing \cite{MR58655} initiated the investigation of these families of subsets of $\N$ and discovered that $\Sigma_\star$ had full measure, in other words, that
almost all sets of integers are totally irreducible; the measure here is the natural product probability measure on $\mathcal{P}(\N)$. 

In this article, we will use  a stronger notion of irreducibility, 
which describes a set, every  ``small'' alteration (deletion or insertion of small subsets) of which
is again irreducible.
Following \cite[Definition 1.3]{MR4607615}, given $\mathcal{I}\subseteq \mathcal{P}(\mathbb{N})$,
thought to be a family of small subsets, we denote by $\Sigma(\mathcal{I})$ the family of $\mathcal{I}$\emph{-irreducible} sets, namely, 
$$
\Sigma(\mathcal{I})=\left\{A\subseteq \mathbb{N}: A\neq_{\mathcal{I}} (B+C) \text{ for all }B,C\subseteq \mathbb{N} \text{ with }|B|,|C|\ge 2\right\},
$$
where $A\neq_{\mathcal{I}} (B+C)$ stands for $A\bigtriangleup (B+C)\notin \mathcal{I}$. 
We denote by $\mathcal{R}(\mathcal{I})=\Pn\setminus\Sigma(\I)$ its complement.
In particular, 
$$
\Sigma=\Sigma(\{\emptyset\}) \qquad \text{ and }\qquad \Sigma_\star=\Sigma(\mathrm{Fin})
$$
where $\mathrm{Fin}$ stands for the family of all finite subsets of $\N$.

Since $\mathcal{I}$ represents a family of small sets,
this notion is attractive when $\mathcal{I}$ is an \emph{ideal} on $\mathbb{N}$, that is, a nonempty family of subsets stable under subsets and finite unions.
A family $\mathcal{I}\subseteq\Pn$ is \emph{admissible} 
if it is stable under finite modifications, i.e., 
$A\bigtriangleup F\in \I$ for all $A \in \I$ and $F \in \mathrm{Fin}$. 
An admissible family which is an ideal, equivalently an ideal containing $\mathrm{Fin}$, is called
an admissible ideal. 

Examples of ideals are the trivial family $\{\emptyset\}$, the family of finite sets $\mathrm{Fin}$, the summable ideal $$\mathcal{I}_{1/n}:=\left\{A\subseteq \mathbb{N}: \sum_{n \in A}\frac{1}{n+1}<\infty\right\},$$ the family $\mathcal{Z}$ of asymptotic density zero sets, or the families 
$$
\mathcal{Z}_\alpha=\{A\subseteq \mathbb{N}: |A(n)|=o(n^\alpha) \text{ as }n\to \infty\}
$$
for each $\alpha \in (0,1]$, where $A(n):=A\cap [0,n]$ (of course, $\mathcal{Z}_1=\mathcal{Z}$).

Observe that, if $\mathcal{I}$ is an ideal on $\N$, the binary relation between sets of nonnegative integers 
$
A=_\mathcal{I}B\Longleftrightarrow A\bigtriangleup B\in \mathcal{I}
$ 
is an equivalence relation, and $\Sigma(\mathcal{I})$
is the set of all sets of nonnegative integers whose equivalence class
contains no nontrivial sumsets.  
In addition, $\Sigma(\cdot)$ is nonincreasing on $\mathcal{P}(\Pn)$, i.e., $\Sigma(\mathcal{I})\subseteq \Sigma(\mathcal{J})$ for all families $\mathcal{I}, \mathcal{J}\subseteq \Pn$ with $\mathcal{J}\subseteq \mathcal{I}$. Further $\Sigma(\emptyset)=\Pn$ and $\Sigma(\Pn)=\emptyset$,
so it suffices to study $\Sigma(\I)$ in the case of nonempty proper families $\I$. 

\subsection{Topological and measure framework}
 
In this paper, we focus on 
topological and measure-theoretic properties of $\Sigma(\I)$. 
To this end, we identify $\Pn$ with $\{0,1\}^\N$. Equipped with
the product topology, it is a compact space by Tychonoff's theorem.
In fact, $(\Pn,\bigtriangleup)$ can be regarded as the abelian compact topological group $((\Z/2\Z)^\N,+)$. 
We endow this group with the completion $\lambda$ of the
Bernoulli product probability measure
$$
\big(\tfrac12\delta_0+\tfrac12\delta_1\big)^{\otimes\mathbb N}.
$$
Its restriction to the Borel $\sigma$-algebra is the unique
Haar probability measure, where translation is the action
of $\bigtriangleup$. 
Further, the map 
$\{0,1\}^\N \to [0,1]$ defined by $x \mapsto \sum_i x_i2^{-i-1}$ 
preserves the measure,  
and induces a bijection $\Pn\setminus\mathrm{Fin}\rightarrow (0,1]$ (by the nonterminating dyadic expansion) and a homeomorphism
from the set of all infinite non cofinite subsets of $\N$ to the set of all non dyadic real numbers in $(0,1)$. Since cofinite subsets and dyadic rationals are countable, measure and Baire category statements do not depend on the chosen framework. 
Recall that $\{0,1\}^\N$ is a zero-dimensional Polish space, that is, a separable and completely metrizable topological space with a base of clopen sets. 

Recall that a subset of a topological space is
\emph{nowhere dense} if its closure has empty interior, \emph{meager} if it is a
countable union of nowhere dense sets, and \emph{comeager} if its complement is
meager. A set has the \emph{Baire property} if it differs from an open set by a
meager set. By the Baire Category Theorem \cite[Theorem~8.4]{MR1321597}, in a
Polish space every comeager set is dense. 
Lastly, given a pointclass $\Gamma\subseteq \mathcal{P}(\mathcal{P}(\mathbb{N}))$ (for instance, $\Gamma=\{S\subseteq \mathcal{P}(\N): S \text{ is a }G_\delta\text{ subset}\}$), we say that $S\subseteq \mathcal{P}(\N)$ is $\Gamma$\emph{-complete} if $S \in \Gamma$ and $\mathcal{P}(\N)\setminus S\notin \Gamma$, cf.  \cite[Definition 22.9 and Theorem 22.10]{MR1321597}.  

\section{Main results}
The following has been shown in \cite[Theorem 1.4]{MR4607615}.
\begin{thm}\label{thm:original}
    $\Sigma(\mathcal{Z}_\alpha)$ is a comeager subset of $\mathcal{P}(\mathbb{N})$ for all $\alpha \in (0,\nicefrac{1}{3})$.
\end{thm}
It was left open in \cite[Section 3]{MR4607615} whether the smaller subset $\Sigma(\mathcal{Z})$ is comeager as well. Below, we provide an affirmative answer to this question, hence generalizing Theorem \ref{thm:original}. 
To this aim, given a positive real sequence $(w_n)_{n\ge 0}$ such that $\sum_n w_n=\infty$ and $\lim_n w_n/\sum_{k\le n}w_k=0$, define the Erd\H{o}s--Ulam ideal 
$$
\mathcal{Z}_w:=\left\{A\subseteq \mathbb{N}: \lim_{n\to \infty}\frac{\sum_{k \in A\cap [0,n]}w_k}{\sum_{k \in [0,n]}w_k}=0\right\},
$$
see e.g. \cite[Section 1.13]{MR1711328} or \cite{MR3391516, MR748847}. 
It is immediate that, if $w_n=1$ for each $n \in \mathbb{N}$, then $\mathcal{Z}_w=\mathcal{Z}$. 
\begin{thm}\label{thm:main1}
Let $(w_n)_{n\ge 0}$ be a nonincreasing positive real sequence with $\sum_n w_n=\infty$. Then $\Sigma(\mathcal{Z}_w)$ is a comeager subset of $\mathcal{P}(\mathbb{N})$.
\end{thm}

Now we turn to $\Sigma(\I)$ for general classes $\I$ (hence, not necessarily ideals).
We obtain a number of basic topological properties.
\begin{prop}\label{prop:complexities}
    Pick a nonempty family $\mathcal{I}\subseteq\Pn$. Then the following hold\textup{:}
     \begin{enumerate}[label={\rm (\roman*)}]
        \item \label{item:1prop} If $\mathcal{I}$ is $F_\sigma$ then $\Sigma(\mathcal{I})$ is $G_\delta$\textup{;}
        \item \label{item:2prop} If $\mathcal{I}$ is analytic then $\Sigma(\mathcal{I})$ is coanalytic\textup{;}
        \item \label{item:3prop} $\mathcal{P}(\mathbb{N})\setminus \Sigma(\mathcal{I})$ is dense\textup{;}
        \item \label{item:4prop} If $\Sigma(\mathcal{I})$ is comeager then $\Sigma(\mathcal{I})$ is not $F_\sigma$\textup{;}
         \item \label{item:5prop} If $\mathcal{I}$ is admissible then $\Sigma(\mathcal{I})$ is either empty or invariant under finite modifications \textup{(}hence dense\textup{)}\textup{;}
         \item \label{item:6prop} If $\mathcal{I}$ is admissible and $\Sigma(\mathcal{I})$ has the Baire property then $\Sigma(\mathcal{I})$ is either meager or comeager\textup{.}
         \item \label{item:7prop} If $\Sigma(\mathcal{I})$ is comeager, then $\I$ is meager.
    \end{enumerate}
\end{prop}

In view of these results, one might conjecture that $\Sigma(\mathcal{I})$ is comeager for all ideals $\mathcal{I}$ with sufficiently low topological complexity. The following result shows that this is false, even for well-behaved $F_\sigma$ ideals. 
Here, an ideal $\mathcal{I}$ is said to be \emph{translation-invariant} if $S+k:=\{n \in \mathbb{N}: n-k \in S\} \in \mathcal{I}$ for all $S \in \mathcal{I}$ and $k \in \mathbb{Z}$. 
\begin{thm}\label{thm:translationinvariant}
There exists an admissible proper translation-invariant $F_\sigma$ ideal
$\mathcal I$ on $\mathbb N$ such that $\Sigma(\mathcal I)=\emptyset$. 
\end{thm}

Alternatively, one might conjecture that if $\mathcal{I}$ is a meager ideal then $\Sigma(\mathcal{I})$ is either meager or comeager (note that, thanks to Proposition \ref{prop:complexities}, this holds if $\mathcal{I}$ is an admissible analytic ideal). Again, this turns out to be false. 
\begin{thm}\label{thm:meagercounterexample}
There exists an admissible meager ideal $\mathcal I$ on $\mathbb N$ such that
$\Sigma(\mathcal I)$ is neither meager nor comeager in
$\mathcal P(\mathbb N)$.
\end{thm}

Next, relying on Proposition \ref{prop:complexities}, we find the exact topological complexity of certain families of $\mathcal{I}$-irreducible sets.  In a related result, Schmerl proved in \cite{MR1635608} that the family $\{(A-A) \cap \mathbb{N}: A\subseteq \mathbb{N}\}$ is analytic but not coanalytic (hence, in particular, not Borel). 
\begin{thm}\label{thm:main2}
Let $\mathcal{I}$ be an admissible nonempty $F_\sigma$ subset of $\mathcal{P}(\mathbb{N})$. 
Then $\Sigma(\mathcal{I})$ is 
either empty or 
dense $G_\delta$-complete \textup{(}hence comeager\textup{)}.
\end{thm}

Note that $\Sigma(\mathcal{I})$ may well be empty even under the hypotheses of Theorem \ref{thm:main2}, see Theorem \ref{thm:translationinvariant}. Moreover, establishing $\Sigma(\mathcal{I})\neq\emptyset$ is not always easy when it is true:
 in fact, an elementary argument shows that any infinite set $A\subseteq \N$ with increasing enumeration $(a_n: n \in \N)$ satisfying $\lim_n(a_{n+1}-a_n)=\infty$ is totally irreducible, i.e., $A \in \Sigma(\Fin)$, but there does not seem
to be a simpler argument than the probabilistic one for the analogous nonemptiness result for, e.g., $\Sigma(\mathcal{I}_{1/n})$.

Combining Proposition \ref{prop:complexities} and Theorem \ref{thm:main1}, we obtain the following. 
\begin{cor}\label{cor:irreducible}
    Each of the families $\Sigma$, $\Sigma_\star$, and $\Sigma(\mathcal{I}_{1/n})$ is dense $G_\delta$-complete \textup{(}hence comeager\textup{)}, and with empty interior. 
\end{cor}
In particular, the family $\mathcal{R}=\Pn\setminus\Sigma$ of all nontrivial sumsets is $F_\sigma$-complete.

In the case of the ideal $\mathcal{Z}$, we get the following. 
\begin{cor}\label{cor:irreducibleZ}
    $\Sigma(\mathcal{Z})$ is comeager, coanalytic, not $F_\sigma$, and with empty interior. 
\end{cor}

It would be interesting to calculate the exact topological complexity of $\Sigma(\mathcal Z)$, hence, in particular, whether it is Borel or not. 
\begin{question}\label{q:BorelSigmaZ}
Is $\Sigma(\mathcal Z)$ Borel? 
\end{question}

Although we do not have an answer to Question \ref{q:BorelSigmaZ}, we conclude by showing, in particular, that the complexity of $\Sigma(\mathcal Z)$ is at least $G_{\delta\sigma}$.  
\begin{thm}\label{thm:notFsigmadelta}
Let $(w_n)_{n\ge 0}$ be a nonincreasing positive real sequence with $\sum_n w_n=\infty$. Then
    $\Sigma(\mathcal Z_w)$ is 
    not $F_{\sigma\delta}$. 
\end{thm}

The proofs of our results are given in the next section. 

\section{Proofs}

First, we need some preliminaries. To this end, define the family 
$$
\mathcal{Z}_{-}
:=
\left\{
A\subseteq\mathbb{N}:
\liminf_{n\to\infty}\frac{|A(n)|}{n}=0
\right\}.
$$
Notice that $\mathcal{Z}_{-}$ need not be an ideal (more precisely, $\mathcal{Z}_{-}$ is stable under subsets but not under finite unions).  
Hereafter, we will use the following slight abuse of notation common in probability theory:
given a predicate $P(X)$ on a set $X\in\Pn$,
we write $\lambda(P(X))$ for $\lambda(\{X\in \Pn : P(X)\})$. Also, we define
$$
\eta:=\frac{1}{100} 
\qquad \text{ and }\qquad 
I_h(n):=\left(\left(1-\frac{1}{2^{h}}\right)n,\,n\right]\cap\N
$$
for all integers $h,n \in \mathbb{N}$. In addition, for each $t \in (0,1)$, set 
$$
H(t):=-t\log_2t-(1-t)\log_2(1-t).
$$
Observe that $\lim_{t\to 0}H(t)=0$. In particular, for each $t \in (0,1)$, it is possible to fix $\varepsilon_t \in (0,\nicefrac{1}{2})$ such that 
$H(\varepsilon_t)<t/8$. 

\begin{lem}\label{lem:template}
The following properties are satisfied by $\lambda$-almost every set $X\subseteq\N$\textup{:}
\begin{enumerate}[label=\textup{(\roman*)}]
\item \label{item:lemma31i} for every finite $F\subseteq\mathbb Z$, every $G\subseteq F$ and every
$h \in \mathbb{N}$, 
\begin{equation}\label{eq:pattern}
\lim_{n\to \infty} \frac{\big|\{t\in I_h(n):\ t+F\subseteq\N \text{ and } X\cap(t+F)=t+G\}\big|
}{|I_h(n)|}=\frac{1}{2^{|F|}}\,\textup{;}
\end{equation}
\item \label{item:lemma31ii} for every integer $h \ge 1$ and for every sufficiently large
$n \in \mathbb{N}$, there are no sets $P,Q\subseteq[0,n]$ satisfying 
\begin{equation}\label{eq:noapprox}
|P|,|Q|\le\varepsilon_{2^{-h}}\,n\qquad\text{and}\qquad
\big|\big(X\bigtriangleup (P+Q)\big)\cap I_h(n)\big|\le\eta\,|I_h(n)|.
\end{equation}
\end{enumerate}
\end{lem}
 \begin{proof}
\ref{item:lemma31i} Fix $G\subseteq F\subseteq \mathbb{Z}$ and $h \in \mathbb{N}$. If $F=\emptyset$ the claim is clear. Otherwise, set $r:=\max F-\min F$. At most $|\min F|$ elements
$t\in I_h(n)$ violate $t+F\subseteq\N$, so we may restrict to the others. For
these, the random variables $$\xi_t=\xi_t(X):=\mathbf 1\{X\cap(t+F)=t+G\}$$ have mean
$2^{-|F|}$.

For $t$ ranging over a fixed residue class modulo $r+1$, the sets
$t+F$ are pairwise disjoint, so inside each of the $r+1$ classes the $\xi_t$
are independent; each class has at most $2n/(r+1)$ elements if $n\ge r+1$.
Combined with the pigeonhole principle,
Hoeffding's inequality \cite[Theorem~1]{Hoeffding1963}, applied in each class with deviation
$n^{3/4}/(r+1)$, shows that 
the numerator in \eqref{eq:pattern} differs 
from $2^{-|F|}|I_h(n)|$ by more than $n^{3/4}+|\min F|$ with probability at
most $$p_n=\frac{2(r+1)}{e^{\sqrt{n}/(r+1)}}.$$ 
Observe that $\sum_n p_n<\infty$. Hence, the conclusion follows by the
Borel--Cantelli lemma and a countable intersection over $(F,G,h)$.
 
\ref{item:lemma31ii} Fix a positive integer $h$ and set $\varepsilon:=\varepsilon_{2^{-h}}$. Observe that
$$
\forall n \in \mathbb{N}, \forall \theta \in [0,\nicefrac{1}{2}], \qquad \sum_{i=0}^{\theta n}\binom ni\le2^{H(\theta)n}.
$$ 
Hence there are at most $2^{2H(\varepsilon)(n+1)}$ pairs
$(P,Q)$ of subsets of $[0,n]$ of size at most $\varepsilon n$. For a fixed
pair, set $R:=(P+Q)\cap I_h(n)$ and note that the coordinates of $X$ on $I_h(n)$ are
independent and each with probability $\nicefrac{1}{2}$, so
$$\lambda\big(|(X\sd R)\cap I_h(n)|\le\eta|I_h(n)|\big)\le \frac{1}{2^{(1-H(\eta))|I_h(n)|}}.$$
Since $|I_h(n)|\ge 2^{-h}n$, $H(\eta)<\nicefrac{1}{10}$, and $2H(\varepsilon)\le 2^{-h-2}$, the
probability that some pair satisfies \eqref{eq:noapprox} is at most
$2^{-n/2^{h+1}+1}$. As in the previous point, the latter sequence is summable in $n$, hence the conclusion follows by the Borel--Cantelli lemma and a
countable intersection over $h\ge 1$.
\end{proof}
The following lemma is not probabilistic but demonstrates some rigidity properties of a set satisfying the generic properties
outlined in Lemma \ref{lem:template}. 
Roughly, a sumset cannot stay essentially inside $X$ below some large $n$ while
coinciding essentially with $X$ on the final segment $I_h(n)$.

\begin{lemma}
\label{lem:rigid}
Pick a set $X\subseteq \mathbb{N}$ satisfying the conditions of Lemma~\ref{lem:template}, and fix an integer $\ell\ge 1$. 
Then there are no sets $B,C\subseteq\N$ with $|B|,|C|\ge2$ such that 
\begin{equation}\label{eq:rigid}
\liminf_{n\to \infty}\frac{\big|((B+C)\setminus X)(n)\big|+\big|(X\sd (B+C))\cap I_\ell(n)\big|}{n}=0.
\end{equation}
\end{lemma}

\begin{proof}
Let us suppose for the sake of contradiction that there exists a strictly increasing sequence $(n_j)$ in $\mathbb{N}$ and sets $B,C\subseteq \mathbb{N}$ with $|B|,|C|\ge 2$ such that 
$$
\big|(S\setminus X)(n_j)\big|=o(n_j)
\qquad\text{and}\qquad
\big|(X\sd S)\cap I_j\big|=o(n_j).
$$
as $j\to \infty$, where $S:=B+C$ and $I_j:=I_\ell(n_j)$ for all $j \in \mathbb{N}$. Define also $c:=2^{-\ell}$ and 
$$
\forall j \in \mathbb{N}, \qquad e_j:=|(S\setminus X)(n_j)|+|(X\sd S)\cap I_j|,
$$
and observe that $\lim_j e_j/n_j=0$. 

\medskip

\emph{$B$ and $C$ are infinite.} Suppose that $B$ is finite and put
$m:=1+\max B$. By \eqref{eq:pattern} with $F=\{-m,\dots,m\}$, $h=\ell$ and $G=\{0\}$,
the interval $I_j$ contains at least $2^{-(2m+1)}cn_j-o(n_j)$ points $t$ with
$X\cap[t-m,t+m]=\{t\}$. If such a $t$ is not in $S$, it lies in
$(X\sd S)\cap I_j$. Otherwise $t=b+y$ with $b\in B$, $y\in C$; choosing
$b'\in B\setminus\{b\}$, the point $t':=b'+y$ lies in $S$, satisfies
$0<|t'-t|<m$, hence $t'\notin X$, i.e., 
$t'\in(S\setminus X)(n_j+m)$. Each $t'$ arises from at most $2m+1$ points $t$.
Therefore
\begin{align*}
2^{-(2m+1)}cn_j-o(n_j)&\le|(X\sd S)\cap I_j|+(2m+1)\big(|(S\setminus X)(n_j)|+m\big)\\
&\le(2m+1)(e_j+m),
\end{align*}
a contradiction. By symmetry, $C$ is infinite as well.

\medskip
 
\emph{$B(n_j)$ and $C(n_j)$ are small.} Fix $k\in\N$ with
$2^{-(k+1)}\le\varepsilon_c/2$, let $b_0<\dots<b_k$ be the first $k+1$ elements
of $B$ and $F:=\{b_i-b_0:\ i\le k\}\subseteq\N$. At most $b_k$ elements
$y\in C(n_j)$ satisfy $y+b_k>n_j$. If $y\in C$ and $y+b_k\le n_j$, the points
$y+b_0,\dots,y+b_k$ belong to $S(n_j)$; if none of them lies in $S\setminus X$,
then $t:=y+b_0\in[0,n_j]$ satisfies $t+F\subseteq X$, and $y\mapsto t$ is
injective. Each point of $S\setminus X$ can be written as $y+b_i$ in at most
$k+1$ ways. By \eqref{eq:pattern} with $G=F$ and $h=0$, we obtain 
\[
|C(n_j)|\le b_k+1+2^{-(k+1)}n_j+o(n_j)+(k+1)\,|(S\setminus X)(n_j)|
\le\varepsilon_c\,n_j
\]
for all large $j$. Exchanging the roles of $B$ and $C$ gives
$|B(n_j)|\le\varepsilon_cn_j$.

\medskip
 
\emph{Conclusion.} 
Put $P:=B(n_j)$ and $Q:=C(n_j)$. Every element of
$S\cap[0,n_j]$ is a sum of elements of $P$ and $Q$, so $S\cap I_j=(P+Q)\cap I_j$
and
\[
\big|\big(X\sd(P+Q)\big)\cap I_j\big|=|(X\sd S)\cap I_j|\le e_j\le\eta\,cn_j\le\eta|I_j|
\]
for all large $j$. Together with $|P|,|Q|\le\varepsilon_cn_j$, this contradicts
Lemma~\ref{lem:template}\ref{item:lemma31ii} with $h=\ell$ and $n=n_j$.
\end{proof}
The following provides a strengthening of \cite[Theorem 1.5]{MR4607615}. 
\begin{prop}\label{prop:lowerdensityrobust}
$\lambda(\Sigma(\mathcal{Z}_-))=1$. 
\end{prop}

\begin{proof}
By Lemma~\ref{lem:template}, it suffices to show that every
set $X$ which satisfies the two properties of that lemma
belongs to $\Sigma(\mathcal{Z}_-)$.
Let $X$ be as in Lemma~\ref{lem:template}, and suppose that
$D:=X\sd S\in\mathcal{Z}_-$ for some $S=B+C$ with $|B|,|C|\ge2$. Choose $n_j\to\infty$
with $|D(n_j)|=o(n_j)$. 
Since $S\setminus X\subseteq D$ and $(X\mathbin{\triangle} S)\cap I_1(n_j)\subseteq D(n_j)$,
we have
\[
\big|(S\setminus X)(n_j)\big|+\big|(X\mathbin{\triangle} S)\cap I_1(n_j)\big|
\le 2\,|D(n_j)|=o(n_j).
\]
Hence the liminf in \eqref{eq:rigid} with $\ell=1$ vanishes, which contradicts
Lemma~\ref{lem:rigid}.
\end{proof}

Before we proceed to the proof of Theorem \ref{thm:main1}, we show another intermediate claim. Hereafter, if $(w_n)_{n\ge 0}$ is a positive real sequence with $\sum_n w_n=\infty$, we write 
\begin{equation}\label{eq:defW}
\forall n \in \mathbb{N}, \qquad 
W(n):=\sum_{k=0}^n w_k.
\end{equation}

\begin{lem}\label{lem:wdecreasing}
    Let $(w_n)_{n\ge 0}$ be a nonincreasing positive real sequence with $\sum_n w_n=\infty$. Pick sequences $(a_m)$ and $(b_m)$ of positive integers such that
\begin{equation}\label{eq:condlemwincreasing}
a_m<b_m,\qquad \lim_{m\to \infty}a_m=\infty,
\qquad\text{and}\qquad
\limsup_{m\to \infty}\frac{W(a_m)}{W(b_m)}<1.
\end{equation}
Then
\begin{equation}\label{eq:weightedclaim}
\forall D \in \mathcal{Z}_w, \qquad 
\lim_{m\to \infty} \min_{a_m\le n\le b_m}\frac{|D(n)|}{n}=0. 
\end{equation}
\end{lem}
\begin{proof}
    Suppose that \eqref{eq:condlemwincreasing} holds, and pick $D \subseteq \mathbb{N}$ such that \eqref{eq:weightedclaim} fails. 
    Then there exist $\varepsilon>0$ and a strictly increasing
sequence $(m_j)$ of positive integers such that 
$
|D(n)|\ge\varepsilon n
$ 
for every $n\in[a_{m_j},b_{m_j}]$ and every $j\in \mathbb{N}$. 
Since $(w_n)$ is nonincreasing,
summation by parts gives
\begin{displaymath}
\begin{split}
\forall j \in \mathbb{N}, \qquad 
\sum_{d\in D(b_{m_j})}w_d
={}&|D(b_{m_j})|w_{b_{m_j}}
+\sum_{t=0}^{b_{m_j}-1}|D(t)|(w_t-w_{t+1})\\
\ge{}&
\varepsilon b_{m_j}w_{b_{m_j}}
+\varepsilon\sum_{t=a_{m_j}}^{b_{m_j}-1}
t(w_t-w_{t+1})\\
={}&
\varepsilon\left(
a_{m_j}w_{a_{m_j}}+\sum_{t=a_{m_j}+1}^{b_{m_j}}w_t
\right)\\
\ge{}&
\varepsilon\bigl(W(b_{m_j})-W(a_{m_j})\bigr).
\end{split}
\end{displaymath}
It follows that 
$$
\forall j \in \mathbb{N}, \qquad 
\frac{\sum_{d\in D(b_{m_j})}w_d}{W(b_{m_j})}
\ge
\varepsilon
\left(
1-\frac{W(a_{m_j})}{W(b_{m_j})}
\right).
$$
By \eqref{eq:condlemwincreasing}, the right-hand side is bounded below by a positive constant for all large $j$. This proves that 
$D\notin\mathcal Z_w$. 
\end{proof}

\begin{proof}[Proof of Theorem \ref{thm:main1}]
Thanks to Proposition \ref{prop:lowerdensityrobust}, it is possible to fix a set
$X\in\Sigma(\mathcal Z_-)$. 
For every integer $k\ge1$, choose $h_k>k^2$ such that
\begin{equation}\label{eq:choicehk}
\frac{W(k^2)}{W(h_k)}\le\frac{1}{k}.
\end{equation}
For every $m\ge1$, define the open set 
$$
U_m:=
\bigcup_{k\ge m}\left\{
A\subseteq\mathbb N: A\cap(k,h_k]=X\cap(k,h_k]\right\}.
$$
In addition, $U_m$ is dense. Indeed, fix a nonempty
basic open set 
$V:=\{A\subseteq\mathbb N:A(N)=F\}$,  
where $N\in\mathbb N$ and $F\subseteq[0,N]$. Choosing
$k>\max\{m,N\}$, we may prescribe the coordinates in $(k,h_k]$ to
agree with those of $X$ without changing the coordinates in $[0,N]$.
Therefore $V\cap U_m\neq\emptyset$.

It follows from the Baire category theorem, see e.g. \cite[Theorem 8.4]{MR1321597}, that
$$
U:=\bigcap_{m\ge1}U_m
$$
is comeager. To complete the proof, we claim that
\begin{equation}\label{eq:GsubsetSigma}
U\subseteq\Sigma(\mathcal Z_w).
\end{equation}

To this aim, fix $A\in U$. There exists a strictly increasing sequence $(k_j)$ such
that
\begin{equation}\label{eq:agreementblocks}
A\cap(k_j,h_{k_j}]=
X\cap(k_j,h_{k_j}]
\end{equation}
for every $j$. Suppose, towards a contradiction, that
$A\notin\Sigma(\mathcal Z_w)$. Then there exist
$B,C\subseteq\mathbb N$, with $|B|,|C|\ge2$, such that
$$
D:=A\mathbin{\bigtriangleup}(B+C)\in\mathcal Z_w.
$$
Now, set $a_j:=k_j^2$ and $b_j:=h_{k_j}$ for each $j$, and note that $\lim_ja_j=\infty$ and $\lim_j W(a_j)/W(b_j)=0$ by \eqref{eq:choicehk}. It follows by Lemma \ref{lem:wdecreasing} that there exist integers
$n_j\in[k_j^2,h_{k_j}]$ such that
\begin{equation}\label{eq:Dsmallweighted}
\lim_{j\to \infty}\frac{|D(n_j)|}{n_j}=0.
\end{equation}

On the other hand, \eqref{eq:agreementblocks} implies that 
$\bigl(A\mathbin{\bigtriangleup}X\bigr)(n_j)
\subseteq[0,k_j]$.  
Since $n_j\ge k_j^2$, we obtain
\begin{equation}\label{eq:AXsmall}
\frac{|(A\mathbin{\bigtriangleup}X)(n_j)|}{n_j}
\le
\frac{k_j+1}{k_j^2}
\longrightarrow0 \quad \text{ as }j\to \infty.
\end{equation}
Putting together \eqref{eq:Dsmallweighted} and \eqref{eq:AXsmall}, we get
\begin{displaymath}
\frac{|(X\mathbin{\bigtriangleup}(B+C))(n_j)|}{n_j}
\le{}
\frac{|(X\mathbin{\bigtriangleup}A)(n_j)|}{n_j}+
\frac{|(A\mathbin{\bigtriangleup}(B+C))(n_j)|}{n_j},
\end{displaymath}
which converges to $0$ as $j\to \infty$. 
Therefore 
$
X\mathbin{\bigtriangleup}(B+C)\in\mathcal Z_-,
$ 
contradicting $X\in\Sigma(\mathcal Z_-)$. This proves
\eqref{eq:GsubsetSigma}. Since $U$ is comeager, it follows that
$\Sigma(\mathcal Z_w)$ is comeager as well.
\end{proof}

\medskip

\begin{proof}
[Proof of Proposition \ref{prop:complexities}]
    \ref{item:1prop} 
We first claim that $\mathcal R:=\mathcal{P}(\mathbb{N})\setminus \Sigma$ is $F_\sigma$. For, consider
the map $\varphi:\mathcal P(\mathbb N)^2\to\mathcal P(\mathbb N)$ defined by 
$$
\forall (B,C) \in \mathcal P(\mathbb N)^2, \quad 
\qquad
\varphi(B,C):=B+C.
$$
The map $\varphi$ is continuous. Indeed, for every $n\in\mathbb N$,
whether $n\in B+C$ depends only on the coordinates of $B$ and $C$ in
$[0,n]$.

Now, for all integers $b_0<b_1$ and $c_0<c_1$, set
$$
K_{b_0,b_1,c_0,c_1}
:=
\left\{
(B,C)\in\mathcal P(\mathbb N)^2: b_0,b_1\in B \text{ and } c_0,c_1\in C
\right\}.
$$
Each $K_{b_0,b_1,c_0,c_1}$ is a compact subset of
$\mathcal P(\mathbb N)^2$. Therefore
$\varphi[K_{b_0,b_1,c_0,c_1}]$ is compact and, hence, closed. Moreover,
$$
\mathcal R
=
\bigcup_{\substack{b_0<b_1\\c_0<c_1}}
\varphi[K_{b_0,b_1,c_0,c_1}].
$$
Thus $\mathcal R$ is $F_\sigma$.

Since both $\mathcal I$ and $\mathcal R$ are $F_\sigma$, there exist closed
sets $E_n,F_m\subseteq\mathcal P(\mathbb N)$ such that $\mathcal I=\bigcup_{n}E_n$ and $\mathcal R=\bigcup_{m}F_m$. 
Consider the continuous map $\Delta:\mathcal P(\mathbb N)^2\longrightarrow\mathcal P(\mathbb N)$ defined by 
$$
\forall (E,S) \in \mathcal P(\mathbb N)^2, \quad 
\qquad
\Delta(E,S):=E\mathbin{\bigtriangleup}S.
$$
By the definition of $\mathcal I$-irreducibility,
\begin{displaymath}
\mathcal P(\mathbb N)\setminus\Sigma(\mathcal I)
=\Delta[\mathcal I\times\mathcal R]=\bigcup_{n,m\in\mathbb N}\Delta[E_n\times F_m].
\end{displaymath}
For every $n,m\in\mathbb N$, the set $E_n\times F_m$ is compact,
because it is a closed subset of the compact space
$\mathcal P(\mathbb N)^2$. Hence its image under $\Delta$ is compact
and, therefore, closed. It follows that 
$\mathcal P(\mathbb N)\setminus\Sigma(\mathcal I)$ 
is $F_\sigma$. Therefore $\Sigma(\mathcal I)$ is $G_\delta$.

    \medskip

    \ref{item:2prop} Consider the continuous map $\Phi:\mathcal P(\mathbb N)^3\to\mathcal P(\mathbb N)$ defined by 
$$
\forall (A,B,C) \in \mathcal P(\mathbb N)^3, \quad 
\Phi(A,B,C):=A\mathbin{\bigtriangleup}(B+C).
$$
Indeed, the value of the $n$-th coordinate of $\Phi(A,B,C)$ depends
only on the coordinates of $A,B$, and $C$ in $[0,n]$. Now set
$$
\mathscr X
:=
\left\{
(A,B,C)\in\mathcal P(\mathbb N)^3:
|B|,|C|\ge2
\text{ and }
\Phi(A,B,C)\in\mathcal I
\right\}.
$$
Since $\mathcal I$ is analytic, its continuous preimage
$\Phi^{-1}[\mathcal I]$ is analytic. Moreover, the conditions
$|B|,|C|\ge2$ define an open subset of
$\mathcal P(\mathbb N)^3$. Since analytic sets are stable under
intersections with Borel sets, it follows that $\mathscr X$ is
analytic.

By definition,
$$
\mathcal P(\mathbb N)\setminus\Sigma(\mathcal I)
=
\left\{
A\subseteq\mathbb N:
\exists B,C\subseteq\mathbb N,\ (A,B,C)\in\mathscr X
\right\}.
$$
Thus $\mathcal P(\mathbb N)\setminus\Sigma(\mathcal I)$ is analytic,
being the projection of an analytic set. Therefore
$\Sigma(\mathcal I)$ is coanalytic.

    \medskip
 
\ref{item:3prop} Fix $E\in\mathcal I$. By the definition of
$\mathcal I$-irreducibility,
$$
\{E\mathbin{\bigtriangleup}S:S\in\mathcal R\}
\subseteq
\mathcal P(\mathbb N)\setminus\Sigma(\mathcal I).
$$
Since $S\mapsto E\mathbin{\bigtriangleup}S$ is a homeomorphism,
it is enough to show that $\mathcal R$ is dense. 
To this end, fix a nonempty basic open set
$$
U:=\{A\subseteq\mathbb N:A(k)=F\},
$$
where $k\in\mathbb N$ and $F\subseteq[0,k]$. Define 
$$
B:=\{0,k+1\}
\quad \text{ and }\quad 
C:=F\cup\{k+1,k+2\}.
$$
Then $|B|,|C|\ge2$ and 
$(B+C)(k)=F$.  
Thus $B+C\in U\cap\mathcal R$.

    \medskip

\ref{item:4prop} Suppose that $\Sigma(\mathcal{I})$ is comeager and, thanks to item \ref{item:3prop} its complement is dense. It follows by the Baire Category
Theorem, see e.g. \cite[Theorem 8.4]{MR1321597}, that $\Sigma(\mathcal{I})$ cannot be $F_\sigma$.  

    \medskip

\ref{item:5prop} 
Suppose that $\Sigma(\mathcal{I})\neq\emptyset$ and fix $A\in\Sigma(\mathcal{I})$. We claim that $A\mathbin{\bigtriangleup}F\in\Sigma(\mathcal{I})$ for every $F\in\mathrm{Fin}$. Indeed, otherwise there exist $F\in\mathrm{Fin}$ and
    $B,C\subseteq\mathbb N$ with $|B|,|C|\ge2$ such that 
    $
    (A\mathbin{\bigtriangleup}F)
    \mathbin{\bigtriangleup}(B+C)\in\mathcal I.
    $ 
    By admissibility, 
    $$
    A\mathbin{\bigtriangleup}(B+C)
    =
    F\mathbin{\bigtriangleup}
    \left((A\mathbin{\bigtriangleup}F)
    \mathbin{\bigtriangleup}(B+C)\right)
    \in\mathcal I,
    $$
    which contradicts $A\in\Sigma(\mathcal I)$. 

    Now, fix a nonempty basic open set 
    $
    U:=\{S\subseteq\mathbb N:S(k)=E\},
    $ 
    where $k\in\mathbb N$ and $E\subseteq[0,k]$. Define 
    $
    F:=(A\cap[0,k])\mathbin{\bigtriangleup}E.
    $ 
    Then $F\in\mathrm{Fin}$ and 
    $
    (A\mathbin{\bigtriangleup}F)(k)=E.
    $ 
    Hence, by the above observation, 
    $
    A\mathbin{\bigtriangleup}F
    \in U\cap\Sigma(\mathcal I).
    $ 
    Thus every nonempty basic open subset of $\mathcal P(\mathbb N)$
    has nonempty intersection with $\Sigma(\mathcal I)$.

    \medskip

    \ref{item:6prop} It follows by item \ref{item:5prop} and 
    \cite[Theorem~8.46]{MR1321597}. 

    \medskip

\ref{item:7prop} Set
$S:=\{0,1\}+\{0,1\}=\{0,1,2\}$.
By definition, 
$$
\{E\mathbin{\bigtriangleup}S:E\in\mathcal I\}
\subseteq
\mathcal P(\mathbb N)\setminus\Sigma(\mathcal I)=\mathcal{R}(\mathcal{I}).
$$ 
If $\Sigma(\mathcal I)$ is comeager, $\mathcal{R}(\mathcal{I})$ is
meager. Since $E\mapsto E\mathbin{\bigtriangleup}S$ is a
homeomorphism, $\mathcal I$ is meager.
\end{proof}

\medskip

\begin{rmk}\label{rmk:Zminus}
The family $\mathcal{Z}_-$ shows that measure and category may disagree for $\Sigma(\cdot)$.
On the one hand, $\lambda(\Sigma(\mathcal{Z}_-))=1$ by Proposition~\ref{prop:lowerdensityrobust}. On the other hand,
$\Sigma(\mathcal{Z}_-)$ is meager. Indeed,
\[
\mathcal{Z}_-=\bigcap_{k\ge1}\ \bigcap_{N\in\mathbb N}\ \bigcup_{n\ge N}
\Big\{A\subseteq\mathbb N:\ |A(n)|<\tfrac nk\Big\}
\]
is a $G_\delta$ subset of $\mathcal P(\mathbb N)$:
each inner set is clopen, and each union over $n\ge N$
is open and dense. In fact, any prescribed finite initial
segment can be followed by a sufficiently long interval
containing no elements of $A$. 
Hence $\mathcal Z_-$ is comeager. 
Moreover, $\mathcal{Z}_-$ is stable under finite modifications, so for every $A\in\mathcal{Z}_-$ we
have $A\mathbin{\bigtriangleup}(\{0,1\}+\{0,1\})=A\mathbin{\bigtriangleup}\{0,1,2\}\in\mathcal{Z}_-$,
that is, $A\notin\Sigma(\mathcal{Z}_-)$. Therefore
$\Sigma(\mathcal{Z}_-)\subseteq\mathcal P(\mathbb N)\setminus\mathcal{Z}_-$, which is meager.
Note that $\Sigma(\mathcal{Z}_-)$ is nevertheless dense, by
Proposition~\ref{prop:complexities}\ref{item:5prop}.

This contrasts with $\Sigma(\mathcal{Z})$, which is both comeager (Theorem~\ref{thm:main1}) and of full
measure. 
\end{rmk}

\medskip

\begin{proof}
[Proof of Theorem \ref{thm:translationinvariant}]
Set $m_j:=2^{2^{j+2}}$ for all $j\in \mathbb{N}$, and define
$$
E:=\bigcup_{j \in \mathbb{N}}E_j, 
\quad \text{ where }\quad 
E_j:=\mathbb{N}\cap [2m_j,2m_j+j]\,\text{ for all }j \in \mathbb{N}.
$$

Observe that $2m_{j-1}+j<m_j$ for all $j\ge 1$, and that for every nonempty finite $F\subseteq \mathbb{Z}$, we have $\bigcap_{t\in F}(E+t)$ is infinite, where $E+t:=\{n \in \mathbb{N}: n-t \in E\}$. Indeed, if
$j\ge\max F-\min F$ is sufficiently large, then
\[
\bigcap_{t\in F}(E+t)\supseteq 
\bigcap_{t\in F}(E_j+t)
\supseteq 
[2m_j+\max F,2m_j+j+\min F]\cap \mathbb{N}\neq \emptyset. 
\]
This implies that no finite union of translates of $E^c:=\mathbb{N}\setminus E$ can be cofinite. 

At this point, define the ideal 
$$
\mathcal{I}:=\left\{
D\subseteq\mathbb N:
D\setminus \bigcup_{t\in F}(E^c+t) \in \mathrm{Fin}
\text{ for some nonempty finite }F\subseteq\mathbb Z
\right\}.
$$
Then it is routine to check that $\mathcal{I}$ is a proper admissible translation-invariant $F_\sigma$ ideal.

To complete the proof, it is enough to show that every subset of $\mathbb N$ is not $\mathcal I$-irreducible. 
To this aim, fix
$A\subseteq\mathbb N$ and define $H_j
:=
\{h\in[0,j]:2m_j+h\in A\}$ for all $j \in \mathbb{N}$. Then, set 
$$
B:=\{m_j:j\in\mathbb N\} 
\quad \text{ and }\quad 
C:=\{0,1\}\cup
\bigcup_{j\in\mathbb N}(H_j+m_j).
$$
Clearly, $|B|,|C|\ge2$. 

At this point, we claim that
\begin{equation}\label{eq:agreementE}
(B+C)\cap E=A\cap E.
\end{equation}
Fix $j\in\mathbb N$. The sums of the type $m_j+(m_j+h)=2m_j+h$ with $h \in H_j$ 
produce precisely the elements of $A\cap E_j$. 
It remains to check that no other sum enters $E_j$. If $i\neq k$ and
$r:=\max\{i,k\}$, then $2m_{r-1}+r-1<m_i+m_k+h<2m_r$ whenever $0\le h\le k$. 
Thus every sum
$m_i+(m_k+h)$ with $i\neq k$ lies strictly between two consecutive intervals
forming $E$. Similarly, each of the sums $m_i$ and $m_i+1$, arising
from the two distinguished elements $0,1\in C$, lies outside $E$.
This proves \eqref{eq:agreementE}. 
By the definition of $\mathcal{I}$, we conclude that $A\mathbin{\bigtriangleup}(B+C) \subseteq\mathbb N\setminus E \in \mathcal{I}$. Therefore $\Sigma(\mathcal{I})=\emptyset$. 
\end{proof}

\medskip

\begin{proof}
[Proof of Theorem \ref{thm:meagercounterexample}]
Let $\mathcal J$ be a maximal ideal on $\mathbb N$ such that $\mathrm{Fin}\subseteq \mathcal{J}$ (that is, the complement of a free ultrafilter). Define $m_0:=0$, $m_n:=2^{2^n}$ for $n\ge1$, and $I_n:=[m_n,m_{n+1})$. For every $A\subseteq\mathbb N$ and $n \in \mathbb{N}$, set $\mu_n(A):=|A\cap I_n|/|I_n|$, 
and define
\begin{equation}\label{eq:auxiliaryideal}
\mathcal K
:=
\left\{
A\subseteq\mathbb N:
\mathcal{J}\text{-}\lim_{n\to \infty}\mu_n(A)=0
\right\}.
\end{equation}
Here, for a real sequence $(x_n)$ and $a\in\mathbb R$,
the notation $\mathcal J\text{-}\lim_{n\to\infty}x_n=a$
means that 
$
\{n\in\mathbb N:|x_n-a|\ge\varepsilon\}\in\mathcal J
$ for all $\varepsilon>0$. 
Since $\mathcal J$ is maximal, every bounded real sequence
has a unique $\mathcal J$-limit. 
Then $\mathcal K$ is a proper ideal containing $\mathrm{Fin}$.

\medskip

First, we claim that 
\begin{equation}\label{eq:ksubsetlower}
\mathcal K\subseteq\mathcal Z_-.
\end{equation}
Indeed, if $A\in\mathcal K$, there exists a strictly increasing
sequence $(n_j)$ such that
\[
\forall j \in \mathbb{N}, \qquad 
\mu_{n_j}(A)\le2^{-j}
\qquad\text{and}\qquad
\frac{m_{n_j}}{m_{n_j+1}}\le 2^{-j}.
\]
It follows that
\[
\frac{|A(m_{n_j+1}-1)|}{m_{n_j+1}-1}
\le
\frac{m_{n_j}}{m_{n_j+1}-1}
+
\frac{|I_{n_j}|}{2^j(m_{n_j+1}-1)}
\longrightarrow 0
\quad \text{ as }j\to \infty, 
\]
which proves \eqref{eq:ksubsetlower}.

\medskip

Second, set
$E:=2\mathbb N$ and $O:=2\mathbb N+1$. We claim that both 
\[
\mathcal K
\quad \text{ and }\quad 
\mathcal K_{\mathrm o}
:=
\mathcal K\cap\mathcal P(O)
\,\,\text{ are not meager},
\]
where the latter is regarded as a subset of
$\mathcal P(O)$. Fix
$Q\in\{\mathbb N,O\}$ and let $(P_k)$ be an arbitrary
partition of $Q$ into finite sets. For each $L\subseteq\mathbb N$, define
\[
\nu(L)
:=
\mathcal{J}\text{-}\lim_{n\to\infty}
\frac{|P_L\cap I_n\cap Q|}{|I_n\cap Q|}, 
\quad \text{ where }
P_L:=\bigcup_{k\in L}P_k.
\]
Then
$\nu$ is a finitely additive probability measure on
$\mathcal P(\mathbb N)$ and 
$
\nu(\{k\})=0 
$ for each $k \in \mathbb{N}$. 
By successively splitting infinite sets into two infinite pieces, we
obtain that there exist decreasing infinite sets $(L_j)$ such that
\begin{equation}\label{eq:nullsets}
\nu(L_j)\le2^{-j}.
\end{equation}
Choose distinct $k_j\in L_j$ and set $L:=\{k_j:j\in\mathbb N\}$.
Since 
$
L\subseteq\{k_0,\ldots,k_{j-1}\}\cup L_j,
$ 
it follows by Inequality \eqref{eq:nullsets} that $\nu(L)=0$. Consequently,
$P_L\in\mathcal K$ if $Q=\mathbb N$, and
$P_L\in\mathcal K_{\mathrm o}$ if $Q=O$, while $P_L$
contains infinitely many members of the given partition. 
By Talagrand's characterization 
\cite[Theorem 21]{Talagrand}, an ideal $\mathcal M$ on
$\mathbb N$ is meager if and only if there exists a partition
$(S_n)$ of $\mathbb N$ into finite intervals such that every
$A\in\mathcal M$ contains at most finitely many of the sets $S_n$.
Hence both $\mathcal K$ and
$\mathcal K_{\mathrm o}$ are not meager.

\medskip

At this point, define 
\begin{equation}\label{eq:meagerideal}
\mathcal I:=\mathcal K\cap\mathcal L, 
\quad \text{ where }\quad 
\mathcal L
:=
\left\{
A\subseteq\mathbb N:
|A\cap2\mathbb N|<\infty
\right\}.
\end{equation}
The family $\mathcal L$ is an $F_\sigma$ ideal, hence meager. It follows that 
$\mathcal I$ is meager as well.

\medskip

Next, we show that $\Sigma(\mathcal I)$ is not meager. By Proposition
\ref{prop:lowerdensityrobust}, it is possible to pick 
$X\in\Sigma(\mathcal Z_-)$. In view of
\eqref{eq:ksubsetlower}, we have $X\in\Sigma(\mathcal K)$. Moreover,
\begin{equation}\label{eq:translatek}
X\mathbin{\bigtriangleup}\mathcal K
:=
\left\{
X\mathbin{\bigtriangleup}Y:Y\in\mathcal K
\right\}
\subseteq\Sigma(\mathcal K).
\end{equation}
Indeed, if $Y\in\mathcal K$ and
\[
(X\mathbin{\bigtriangleup}Y)
 \mathbin{\bigtriangleup}(B+C)\in\mathcal K,
\]
then $X\mathbin{\bigtriangleup}(B+C)\in\mathcal K$, which is a
contradiction. 
Since $Y \mapsto X \mathbin{\bigtriangleup} Y$ is a homeomorphism for each $X$, and since 
$\mathcal K$ is not meager, \eqref{eq:translatek} implies that
$\Sigma(\mathcal K)$ is not meager. Finally,
$\mathcal I\subseteq\mathcal K$ gives 
$
\Sigma(\mathcal K)\subseteq\Sigma(\mathcal I), 
$ 
and hence $\Sigma(\mathcal I)$ is not meager.

\medskip

It remains to prove that $\Sigma(\mathcal I)$ is not comeager. Define $\Phi:\mathcal P(\mathbb N)\to\mathcal P(O)$ by 
\[
\forall A \in  \mathcal P(\mathbb N), 
\qquad
\Phi(A):=
A\mathbin{\bigtriangleup}
\bigl(\{0,1\}+(A\cap E)\bigr).
\]
Notice that $\Phi(A)\subseteq O$ for every $A$. Moreover, the map $\mathcal P(\mathbb N)\to \mathcal P(E)\times\mathcal P(O)$ defined by 
$
A\mapsto\bigl(A\cap E,\Phi(A)\bigr)
$ 
is a homeomorphism. 
Since $\mathcal K_{\mathrm o}$ is not meager, it follows that
\begin{equation}\label{eq:nonmeagerpreimage}
\Phi^{-1}[\mathcal K_{\mathrm o}]
\quad\text{is not meager}.
\end{equation}

The set 
$
G:=\{A\subseteq\mathbb N:|A\cap E|\ge2\}
$ 
is open and dense. Thus
$G\cap\Phi^{-1}[\mathcal K_{\mathrm o}]$ is not meager by
\eqref{eq:nonmeagerpreimage}. If $A$ belongs to this intersection,
put
\[
B:=\{0,1\}
\qquad\text{and}\qquad
C:=A\cap E.
\]
Then $|B|,|C|\ge2$ and
\[
A\mathbin{\bigtriangleup}(B+C)
=
\Phi(A)\in\mathcal K_{\mathrm o}
\subseteq\mathcal K\cap\mathcal L
=
\mathcal I.
\]
Consequently,
\[
G\cap\Phi^{-1}[\mathcal K_{\mathrm o}]
\subseteq
\mathcal P(\mathbb N)\setminus\Sigma(\mathcal I).
\]
The complement of $\Sigma(\mathcal I)$ is therefore not meager, so
$\Sigma(\mathcal I)$ is not comeager.
\end{proof}

\medskip

\begin{proof}[Proof of Theorem \ref{thm:main2}]
It follows by items \ref{item:1prop}, \ref{item:3prop}, \ref{item:4prop}, and \ref{item:5prop} in Proposition \ref{prop:complexities}.  
\end{proof}

\medskip

\begin{proof}
    [Proof of Corollary \ref{cor:irreducible}] 
    Recall that $\{\emptyset\}$, $\mathrm{Fin}$, and $\mathcal{I}_{1/n}$ are $F_\sigma$ ideals, with $\{\emptyset\}\subseteq \mathrm{Fin}\subseteq \mathcal{I}_{1/n}\subseteq \mathcal{Z}$. Hence $\Sigma(\mathcal{Z})\subseteq \Sigma(\mathcal{I}_{1/n})\subseteq \Sigma_\star \subseteq \Sigma$. Since $\Sigma(\mathcal{Z})$ is comeager by Theorem \ref{thm:main1}, these families are comeager. 
    The remaining claims follow by items \ref{item:1prop}, \ref{item:3prop}, and \ref{item:4prop} in Proposition \ref{prop:complexities}.  
\end{proof}

\medskip

\begin{proof}
[Proof of Corollary \ref{cor:irreducibleZ}]
Recall that $\mathcal Z$ is an $F_{\sigma\delta}$ subset
of $\mathcal P(\mathbb N)$, hence analytic.
The conclusion follows by Theorem \ref{thm:main1} and
items \ref{item:2prop}, \ref{item:3prop}, and
\ref{item:4prop} in Proposition \ref{prop:complexities}.
\end{proof}

\medskip

\begin{proof}
[Proof of Theorem \ref{thm:notFsigmadelta}]
Consider the set
$$
S_3
:=
\left\{
x\in2^{\mathbb N\times\mathbb N}:
\exists i\in\mathbb N\ \forall N\in\mathbb N\
\exists n\ge N,\ x(i,n)=1
\right\}.
$$
Thus $x\in S_3$ if and only if at least one row of $x$ contains
infinitely many $1$'s. It is well known that $S_3$ is
$G_{\delta\sigma}$-complete (that is, $G_{\delta\sigma}$ and not $F_{\sigma\delta}$), cf. 
\cite[p.~179]{MR1321597}.

For every
$x\in2^{\mathbb N\times\mathbb N}$ and $n\ge 1$, put
$$
\delta_n(x)
:=
\max\left(
\left\{
2^{-(i+1)}:i\le n\text{ and }x(i,n)=1
\right\}
\cup\{0\}
\right).
$$
Then
\begin{equation}\label{eq:deltacharacterization}
x\in S_3
\quad\Longleftrightarrow\quad
\limsup_{n\to\infty}\delta_n(x)>0.
\end{equation}
Indeed, the forward implication is immediate. Conversely, suppose that
$\limsup_n\delta_n(x)>0$. By the pigeonhole principle, there exists a positive integer $r$
such that 
$
\delta_n(x)\ge2^{-r}
$ 
for infinitely many $n$. For every such $n$, there exists $i<r$ such
that $x(i,n)=1$. Hence one of the first
$r$ rows contains infinitely many $1$'s. 
Thus $x\in S_3$, proving
\eqref{eq:deltacharacterization}.

\medskip

Fix a set $X$ as in Lemma \ref{lem:template}, and choose integers
$1\le h_0<h_1<\cdots$ such that
\begin{equation}\label{eq:hnweighted}
\forall n\ge1,\qquad
\frac{W(h_{n-1})}{W(h_n)}\le2^{-n},
\end{equation}
where $W(\cdot)$ is defined in \eqref{eq:defW}. 
For every $x\in2^{\mathbb N\times\mathbb N}$ and $n\ge1$, define
$$
a_n(x):=\min\{a\in[0,h_n]\cap\mathbb N:
W(a)\ge(1-\delta_n(x))W(h_n)\},
$$
and put
$$
J_n(x):=(a_n(x),h_n]\cap\mathbb N
\qquad\text{and}\qquad
A_x:=X\cap\bigcup_{n\ge1}J_n(x).
$$
By the minimality of $a_n(x)$ and the fact that $w_k\le w_0$,
\begin{equation}\label{eq:anweighted}
(1-\delta_n(x))W(h_n)
\le W(a_n(x))
\le(1-\delta_n(x))W(h_n)+w_0.
\end{equation}
Since $0\le\delta_n(x)\le\nicefrac12$, Equations
\eqref{eq:hnweighted} and \eqref{eq:anweighted} imply that
$a_n(x)\ge h_{n-1}$. Hence $J_n(x)\subseteq(h_{n-1},h_n]$, and the
intervals $J_n(x)$ are pairwise disjoint. Moreover, $J_n(x)=\emptyset$
if $\delta_n(x)=0$.

The value of $a_n(x)$ depends only on the finitely many coordinates
$x(0,n),\ldots,x(n,n)$. For every fixed $m\in\mathbb N$, there is at
most one $n\ge1$ for which $m\in(h_{n-1},h_n]$. Thus, the map
$$
\Psi:2^{\mathbb N\times\mathbb N}\to\mathcal P(\mathbb N):x\mapsto A_x
$$
is continuous.

\medskip

Next, we claim that
\begin{equation}\label{eq:reductionclaim}
x\in S_3
\quad\Longleftrightarrow\quad
A_x\in\Sigma(\mathcal Z_w).
\end{equation}
Suppose first that $x\notin S_3$. By
\eqref{eq:deltacharacterization} we have $\lim_n\delta_n(x)=0$.
For every $n\ge1$, Equation \eqref{eq:anweighted} implies
$$
\sum_{k\in A_x(h_n)}w_k
\le W(h_{n-1})+\delta_n(x)W(h_n).
$$
Fix $n\ge2$ and $N\in(h_{n-1},h_n]$.
If $N\le a_n(x)$, then $A_x(N)=A_x(h_{n-1})$, and hence
$$
\frac{\sum_{k\in A_x(N)}w_k}{W(N)}
\le
\frac{W(h_{n-2})+\delta_{n-1}(x)W(h_{n-1})}{W(h_{n-1})}
\le2^{1-n}+\delta_{n-1}(x).
$$
On the other hand, if $N>a_n(x)$, then
$W(N)\ge W(a_n(x))\ge W(h_n)/2$. Hence 
$$
\frac{\sum_{k\in A_x(N)}w_k}{W(N)}
\le
2\left(\frac{W(h_{n-1})}{W(h_n)}+\delta_n(x)\right)
\le2^{1-n}+2\delta_n(x).
$$
It follows that $A_x\in\mathcal Z_w$, and so
$A_x\mathbin{\bigtriangleup}(\{0,1\}+\{0,1\})\in\mathcal Z_w$ since $\Fin\subseteq\mathcal{Z}_w$.
Therefore $A_x\notin\Sigma(\mathcal Z_w)$.

\medskip

Suppose now that $x\in S_3$. By
\eqref{eq:deltacharacterization}, there exist a positive integer $j$
and an infinite set $L\subseteq\mathbb N$ such that
\begin{equation}\label{eq:fixedwindow}
\forall n\in L,\qquad
\delta_n(x)=2^{-j}=:\delta.
\end{equation}
Set $c:=\delta/4=2^{-(j+2)}$, and, for every $n\in L$, put
$a_n:=a_n(x)$ and $a_n':=\lceil a_n/(1-c)\rceil$.
Since $(w_k)$ is nonincreasing and $c\le\nicefrac18$, we have
$$
W(a_n')-W(a_n)
\le(a_n'-a_n)w_{a_n}
\le(2ca_n+1)w_{a_n}
\le2cW(a_n)+w_0.
$$
Together with \eqref{eq:anweighted} and
$(1+\delta/2)(1-\delta)\le1-\delta/2$, this implies, for all sufficiently
large $n\in L$,
\begin{equation}\label{eq:weightedgap}
W(a_n')
\le(1-\delta/2)W(h_n)+3w_0
\le(1-\delta/4)W(h_n).
\end{equation}
In particular, $a_n'<h_n$ for all such $n$.

Suppose, towards a contradiction, that there are $B,C\subseteq\N$ with
$|B|,|C|\ge2$ such that $D:=A_x\sd S\in\mathcal{Z}_w$, where $S:=B+C$. We apply
Lemma~\ref{lem:wdecreasing} to the pairs $a_n'<h_n$, $n\in L$: indeed
$a_n'\ge a_n\ge h_{n-1}\to\infty$ and
$\limsup_{n\in L}W(a_n')/W(h_n)\le1-\frac\delta4<1$. Hence there are integers
$n^*=n^*(n)\in[a_n',h_n]$ and reals $\tau_n\to0$ such that
\begin{equation}\label{eq:good}
\forall n \in L, \qquad 
|D(n^*)|\le\tau_n\,n^*.
\end{equation}
Put $I:=((1-c)n^*,n^*]$. Since $(1-c)n^*\ge(1-c)a_n'\ge a_n$ and
$n^*\le h_n$, we have $I\subseteq J_n(x)$, hence $A_x\cap I=X\cap I$. Recall
also that $A_x\subseteq X$.

Since $A_x\subseteq X$, we have $S\setminus X\subseteq S\setminus A_x\subseteq D$;
and since $A_x\cap I=X\cap I$, we have $(X\mathbin{\triangle} S)\cap I=D\cap I$.
Note that $I=I_{j+2}(n^*)$, because $c=2^{-(j+2)}$. Hence
\[
\big|(S\setminus X)(n^*)\big|+\big|(X\mathbin{\triangle} S)\cap I_{j+2}(n^*)\big|
\le 2\,|D(n^*)|\le 2\tau_n n^*
\qquad(n\in L).
\]
Since $\tau_n\to0$ and $n^*=n^*(n)\ge h_{n-1}\to\infty$ as $n\to\infty$ in $L$, the
liminf in \eqref{eq:rigid} with $\ell:=j+2$ vanishes, which contradicts
Lemma~\ref{lem:rigid}. Therefore $A_x\in\Sigma(\mathcal Z_w)$, proving
\eqref{eq:reductionclaim}.

\medskip

Finally, suppose that $\Sigma(\mathcal Z_w)$ were $F_{\sigma\delta}$.
Since $\Psi$ is continuous, \eqref{eq:reductionclaim} would imply that
$$
S_3
=\left\{x\in2^{\mathbb N\times\mathbb N}:
A_x\in\Sigma(\mathcal Z_w)\right\}
=\Psi^{-1}\left[\Sigma(\mathcal Z_w)\right]
$$
is $F_{\sigma\delta}$ as well, which provides the desired contradiction.
Therefore $\Sigma(\mathcal Z_w)$ is not $F_{\sigma\delta}$.
\end{proof}

\begin{rmk}\label{rmk:weightedhardness}
The proof of Theorem \ref{thm:notFsigmadelta} actually
shows that $\Sigma(\mathcal Z_w)$ is
$G_{\delta\sigma}$-hard for every sequence $(w_n)$
satisfying its hypotheses. This means that, for every
zero-dimensional Polish space $Y$ and every
$G_{\delta\sigma}$ subset $E\subseteq Y$, there exists
a continuous map $f:Y\to\mathcal P(\mathbb N)$ such that
$$
E=f^{-1}\left[\Sigma(\mathcal Z_w)\right].
$$
Indeed, since $S_3$ is $G_{\delta\sigma}$-complete,
there exists a continuous map
$g:Y\to2^{\mathbb N\times\mathbb N}$ such that
$E=g^{-1}[S_3]$. The continuous map $\Psi$ constructed
in the proof satisfies
$S_3=\Psi^{-1}[\Sigma(\mathcal Z_w)]$.
Thus $f:=\Psi\circ g$ has the required property.
\end{rmk}
\section*{Acknowledgment} The first author was funded by the grant 10.55776/PAT4719224 of the Austrian Science Fund (FWF).
\bibliographystyle{amsplain}
\bibliography{sumsets}

\end{document}